\documentclass{article}
\usepackage{fullpage}

\usepackage{amsmath}
\usepackage{amsfonts}
\usepackage{amsthm}
\usepackage{amssymb}
\usepackage{multirow}
\usepackage{algorithm,algorithmic}
\usepackage{mathrsfs}
\usepackage{color,xcolor}
\usepackage{ifpdf}
\usepackage{psfrag}
\usepackage{graphicx,graphics,epsfig,subfigure}
\usepackage{url}
\usepackage{hyperref}
\usepackage{cleveref}
\usepackage{ifpdf}
\usepackage{psfrag}
\usepackage{pgf,tikz}
\usetikzlibrary{arrows.meta, bending, patterns}
\usepackage{url}
\usepackage{hyperref}
\usepackage{xspace}
\usepackage{xargs}
\usepackage{float}
 \usepackage{ booktabs,array}
\usepackage{extarrows}
\usepackage{appendix}
\usepackage{chngcntr}
\AtBeginEnvironment{subappendices}{%
\chapter*{Appendix}
\addcontentsline{toc}{chapter}{Appendices}
\counterwithin{figure}{section}
\counterwithin{table}{section}
}
\usepackage{bm}
\usepackage{dutchcal}
\usepackage{stmaryrd}
\usepackage{comment}

\usepackage[colorinlistoftodos,prependcaption,textsize=tiny]{todonotes}

\newcommand{\remi}[1]{\textcolor{red}{ #1}}

\newtheorem{theorem}{Theorem}[section]

\newtheorem{remark}[theorem]{Remark}

\renewcommand{\div}{\text{ div }}

\newcommand{\R}{\mathbb R}

\renewcommand{\P}{\mathbb P}

\newcommand{\M}{\mathbb{M}}

\newcommand{\EE}{\mathcal E}

\newcommand{\bbf}{\mathbf{f}}

\newcommand{\bbn}{\mathbf{n}}

\newcommand{\bbp}{\mathbf{p}}

\newcommand{\bbu}{\mathbf{u}}

\newcommand{\bbx}{\mathbf{x}}

\newcommand{\bbA}{\mathbf{A}}

\newcommand{\bbJ}{\mathbf{J}}

\newcommand{\bbN}{\mathbf{N}}

\newcommand{\bbU}{\mathbf{U}}

\newcommand{\bbX}{\mathbf{X}}

\newcommand{\bbPhi}{\bm \Phi}
\newcommand{\bbPsi}{\bm \Psi}

\newcommand{\Id}{\mathbb{I}}

\newcommand{\hbbf}{\hat{\mathbf{f}}}

\newcommand{\dt}{\Delta t}

\usepackage{amsopn}

\newcommand*\xbar[1]{%
  \hbox{%
    \vbox{%
      \hrule height 0.5pt % The actual bar
      \kern0.4ex%         % Distance between bar and symbol
      \hbox{%
        \kern-0.05em%      % Shortening on the left side
        \ensuremath{#1}%
        \kern-0.00em%      % Shortening on the right side
      }%
    }%
  }%
}
\newcommand{\balpha}{{\bm {\alpha}}}

\newcommand{\bsigma}{{\bm {\sigma}}}

\newcommand{\bomega}{{\bm {\omega}}}

\newcommand{\bPsi}{{\bm \Psi}}

\newcommand{\dg}{dG\xspace}
\newcommand{\pampa}{PAMPA\xspace}
\newcommand{\dofs}{\text{DoFs}\xspace}
\newcommand{\dof}{\text{DoF}\xspace}

\usepackage[style=numeric,backend=biber,maxnames=999,isbn=false,url=false, doi=false]{biblatex}
\DeclareSourcemap{
  \maps[datatype=bibtex]{
    \map{
      \step[fieldset=language, null]
    }
  }
}
\begin{document}
%%%%%%%%%%%%%%%%%%%%%%%%%%%%%%%%%%%%%%%%%%%%%%%%%%%%%%%%%%%%%%%%%%%%
	\title{Active Flux and Polygons: why a stabilisation and how to implement boundary conditions.}

\author{R\'emi Abgrall\footnote{This paper is dedicated to Philip L. Roe, for he has so much inspired my research  during whole my career. I hope you would have accepted this short note.}~ and Yongle Liu \\
  Institute of Mathematics, University of Zurich
 }
\date{\today} 
\maketitle
\begin{abstract}
In a recent paper, the authors developed an Active Flux-type algorithm for polygonal meshes (see \cite{BP_Pampa_VEM}). The main difference between Active Flux methods on polygonal meshes and those on Cartesian or simplex meshes is that the basis functions are not explicitly known. These functions describe the solution using point-value approximations on the polygon boundaries and the average value within the polygons. In \cite{BP_Pampa_VEM}, it was shown that directly transposing the algorithm from simplexes to polygons leads to a mild instability, and a remedy was proposed without a rigorous explanation. In the present paper, we clarify this issue and demonstrate that the solution presented in \cite{BP_Pampa_VEM} is indeed correct. Furthermore, another challenge—the implementation of boundary conditions—was previously solved heuristically. Building on \cite{pampa_dg}, we show how to properly implement these boundary conditions. Finally, the method is illustrated through several examples with hybrid meshes.
\end{abstract}
%\tableofcontents
\section{Introduction}
We developed in \cite{BP_Pampa_VEM} a scheme for polygonal meshes that follows the Active Flux methodology \cite{AF1,AF2,AF3,EymannRoe2013,Calhoun2023-ms,AF4,H1,H2,BPPAMPA_1D,Abgrall2025b}. In each polygon, the solution is described by its average value and its point values on the boundary. The vertices of the polygons are included among these interpolation points, while the others correspond to Gauss--Lobatto points on each face.
Compared to \cite{AbgrallLinLiu2025,OEPampa}, which were developed for triangles, and even compared to schemes designed purely for Cartesian meshes (such as \cite{barsukow2025generalizedactivefluxmethod,Barsukow,Duan_AF, AF4, H1,H2}), this scheme requires a rather unconventional filtering term. The reason is as follows: unlike standard elements such as simplexes, there is no reference polygon. Consequently, the user must either explicitly define (and store) the basis functions needed to evaluate the functional approximation for each element, or adopt an alternative approach. In \cite{BP_Pampa_VEM}, we utilize the Virtual Element Method (VEM) framework (see \cite{Acta} for a review and Section \ref{sec:setting} for a brief summary). This VEM space is constructed such that the set of polynomials of degree \(k\) (here \(k=2\)) is a subspace; that is, the VEM space \(V(K)\) contains \(\mathbb{P}_{k}\). This implies that we forgot explicit knowledge of the basis functions, knowing simply that they exist. Instead, we compute the \(L^{2}\) projection of any function from the VEM space onto polynomials of degree \(k\), which is sufficient for approximation purposes. In the scheme of \cite{BP_Pampa_VEM}, we need to compute the gradient of the solution—specifically, the gradient of its projection onto \(\mathbb{P}_{k}\). In doing so, information is lost because \(\mathbb{P}_{k}\) is a proper subspace of the VEM space. This issue does not arise in \cite{AbgrallLinLiu2025} or in other Active Flux-type methods. Therefore, a filtering term is required to damp out the unresolved modes without compromising accuracy. In \cite{BP_Pampa_VEM}, a heuristic solution was found based on pragmatic analysis, though it lacked clear theoretical justification.

In that same paper, we also proposed a method to approximate boundary conditions. While the boundary effect on the average value is well understood, the treatment of the point-value degrees of freedom remains problematic. Although the method proposed in \cite{BP_Pampa_VEM} performs well, its derivation is highly heuristic and not entirely satisfactory, at least from a theoretical standpoint.

The purpose of this note is to justify the choice made in our previous work (in fact, we propose an improved formulation). We also demonstrate how to properly approximate the boundary conditions. This justification is achieved by reinterpreting the \pampa scheme within the Discontinuous Galerkin (\dg) framework, following the approach in \cite{OEPampa}.
The outline of this note is as follows. First, we review the approximation framework and rewrite the \pampa scheme using the \dg formalism combined with projection. This allows us to define the appropriate stabilization term as well as the necessary modifications at the boundaries.

\section{Meshes and VEM approximation setting}\label{sec:setting}
A domain \(\Omega \subset \mathbb{R}^d\) is covered by a mesh of non-overlapping polygons \(\{K\}\). The only assumption made on each polygon \(K\) is that it is star-shaped with respect to a point denoted by \(\mathbf{x}^\star_K\) (or simply \(\mathbf{x}^{\star }\) when the context is clear). For the sake of simplicity, we assume \(d=2\) here. The family of polygons is denoted by \(\mathcal{T}=\{K\}\). The set of faces, denoted by \(\mathcal{E}\), consists of the non-empty intersections of polygons. In other words, if \(K, K' \in \mathcal{T}\) and \(e = K \cap K' \neq \emptyset\), then \(e \in \mathcal{E}\). Conversely, for any \(e \in \mathcal{E}\), there exist two polygons \(K\) and \(K^{\prime }\) such that \(e = K \cap K'\). Thus, there is a bijection between the edges and the pairs of polygons with a non-empty intersection. For any \(K\), we denote the set of its faces by \(\mathcal{E}_{K}\) and the set of its vertices by \(\mathcal{V}(K)\). Furthermore, we denote the diameter of \(K\) by \(h_{K}\) and its inner diameter by \(\rho _{K}\). Finally, we assume that the families \(\mathcal{T}\) of polygons under consideration are shape-regular, meaning that the ratio \(h_K / \rho_K\) is bounded independently of \(\mathcal{T}\).

%%
%This being said, for any $K\in \TT$, the set $V_k(K)$ is defined by: for $k\geq 2$
%$$V_k(K)=\{ \bbu\in L^2(K), \forall e\in \EE_K, \bbu_{\mid e}\in \P_k(e) \text{ and }\Delta \bbu\in \P_{k-2}(K)\}.$$
%An element of $V_k(K)$ is defined by the following degrees of freedom:
%\begin{itemize}
%\item The value of $\bbu$ for all the vertices of $K$,
%\item The value of $\bbu$ and the $k-2$ Gauss-Lobatto points of each $e\in \EE_K$,
%\item The moments $m_\balpha(\bbu)$ defined by
%$$m_\balpha(\bbu)=\frac{\int_K \mathcal{M}_\balpha(\bbx)\bbu(\bbx)\; d\bbx}{\vert K\vert}$$
%where $\balpha=(\alpha_1,\alpha_2)$ is a multi-index with $0\leq \vert \balpha\vert =\alpha_1+\alpha_2\leq k-2$, 
%$$\mathcal{M}_\balpha=\Pi_{i=1}^d \big (\frac{x_i-x^\star_i}{h_k}\big )^{\alpha_i}:=\big (\frac{\bbx-\bbx^\star}{h_K}\big )^\balpha$$
%with $\bbx=(x_1,x_2)$.
%\end{itemize}
%In the following, the point values degrees of freedom will be denoted by $\sigma$ where we make the formal identification between the linear form $\bbu\mapsto \bbu({\bbx}_\sigma):=\bbu_\sigma$  and the point itself. It is clear that $\P_k(K)\subset V_k(K)$ and that the inclusion is proper.  From the definition of $V_k(K)$ and the divergence formula, it is possible to define a projection, denoted by $\pi^\nabla$, from $V_k(K)$ onto $\P_k(K)$ by:
%$$\int_K\pi^\nabla(\bbu)\; d\bbx=\int_K \bbu\; d\bbx \text{ and }\int_K {\nabla}v\cdot\big ( \nabla \bbu-\nabla \pi^\nabla \bbu)\; d\bbx=0,\quad \forall v\in \P_k(K), $$
%%

With this in mind, for any \(K \in \mathcal{T}\) and \(k \geq 2\), the space \(V_k(K)\) is defined by
$$V_{k}(K)=\{\mathbf{u}\in L^{2}(K)\,:\,\forall e\in \mathcal{E}_{K},\,\mathbf{u}_{\mid e}\in \mathbb{P}_{k}(e)\text{\ and\ }\Delta \mathbf{u}\in \mathbb{P}_{k-2}(K)\}.$$
An element in \(V_k(K)\) is uniquely determined by the following degrees of freedom:
\begin{itemize}
\item The values of \(\mathbf{u}\) at all vertices of \(K\),
\item The values of \(\mathbf{u}\) at the \(k-1\) internal Gauss--Lobatto points of each edge \(e \in \mathcal{E}_K\),
\item The moments \(m_{{\balpha}}(\mathbf{u})\) defined by
$$m_{{\balpha }}(\mathbf{u})=\frac{1}{|{}K|{}}\int _{K}\mathcal{M}_{{\balpha }}(\mathbf{x})\mathbf{u}(\mathbf{x})\;d\mathbf{x},$$
where \({\balpha}=(\alpha_1,\alpha_2)\) is a multi-index satisfying \(0 \leq \vert{}{\balpha}\vert{} = \alpha_1+\alpha_2 \leq k-2\), and
$$\mathcal{M}_{{\balpha }}(\mathbf{x})=\prod _{i=1}^{d}\left(\frac{x_{i}-x_{i}^{\star }}{h_{K}}\right)^{\alpha _{i}}:=\left(\frac{\mathbf{x}-\mathbf{x}^{\star }}{h_{K}}\right)^{{\balpha }}$$
with \(\mathbf{x}=(x_1,x_2)\).
\end{itemize}
In the following, the point-value degrees of freedom will be denoted by \(\sigma \), where we formally identify the linear form \(\mathbf{u} \mapsto \mathbf{u}(\mathbf{x}_\sigma) := \mathbf{u}_\sigma\) with the point itself. It is clear that \(\mathbb{P}_k(K) \subset V_k(K)\) and that this inclusion is proper. Based on the definition of \(V_k(K)\) and the divergence theorem, one can define an elliptic projection, denoted by \(\pi ^{\nabla }\), from \(V_k(K)\) onto \(\mathbb{P}_k(K)\) via
$$\int _{K}\pi ^{\nabla }(\mathbf{u})\;d\mathbf{x}=\int _{K}\mathbf{u}\;d\mathbf{x}\quad \text{and}\quad \int _{K}\nabla v\cdot \left(\nabla \mathbf{u}-\nabla \pi ^{\nabla }\mathbf{u}\right)\;d\mathbf{x}=0\quad \forall v\in \mathbb{P}_{k}(K).$$
%
%%%%%

The drawback of the space \(V_k(K)\) is that its \(L^{2}\) projection onto \(\P_k(K)\) cannot be computed solely from the knowledge of the degrees of freedom. In contrast, the elliptic projection \(\pi ^{\nabla }\) depends only on these degrees of freedom, as shown by:
\begin{equation*} \left\{ \begin{aligned} & \int_K\pi^{\nabla}\bbu\,d\bbx=\int_K\bbu\,d\bbx,\\ & \int_K\nabla v\cdot\nabla\pi^{\nabla}\bbu\,d\bbx=\int_K \nabla v\cdot\nabla\bbu\,d\bbx=\int_{\partial K}\bbu\nabla v\cdot\bbn\,ds-\int_K\bbu\Delta v\,d\bbx,\quad \forall v\in\P_k(K). \end{aligned}\right. \end{equation*}
To overcome this issue, one must define an alternative space, \(W_k(K)\), which still contains \(\P_k(K)\) but allows for a computable \(L^{2}\) projection. This construction is achieved in two steps. First, we extend \(V_k(K)\) by introducing
$$\tilde{V}_k(K)=\{ \bbu\in L^2(K) \, : \, \forall e\in \EE_K, \, \bbu_{\mid e}\in \P_k(e) \text{ and }\Delta \bbu\in \P_k(K)\},$$
equipped with the same nodal degrees of freedom as \(V_k(K)\), alongside the moments
$$\int_K \bbu \mathcal{M}_\balpha(\bbx)\; d\bbx \quad \text{for } \vert \balpha\vert \leq k.$$
For any $\bbu\in \tilde{V}_k(K)$, its \(L^{2}\) projection can be computed using only its degrees of freedom. Next, letting \(\P_s^{hom}(K)\) denote the set of homogeneous polynomials of degree \(s\), we define
$$W_k(K)=\{ \bbu\in \tilde{V}_k(K) \, : \, \int_K(\bbu-\pi^0\bbu)p_s\; d\bbx=0, \;\forall p_s\in \P_s^{hom}(K), \, s=k-1,k\}.$$
This space shares the exact same degrees of freedom as \(V_k(K)\) and contains \(\P_k(K)\). Ultimately, it can be shown that  \(\dim \widetilde{V_k}(K)=\dim W_k(K)\) and \(\text{dofs}\big(V_k(K)\big)=\text{dofs}\big( W_k(K)\big)\) (see \cite{hitch} for a detailed proof). This is completely transparent to the user.
%%%%%%%%%%%%%%%%%%%%%
\section{\pampa and the \dg formulation}
%%%%%%%%%%%%%%%%%%%%%
\subsection{The \pampa scheme}
For $k=2$, a function $\bbu$ is described as $\{\bbu_\sigma\}$ and the average $\xbar{\bbu}_K$. The nodal degrees of freedom  are the vertices and the mid-points of the edges.
Though we do not use them in practice, and for the discussion only, let us introduce the basis function $\varphi_\sigma$, $\xbar{\varphi}$. A function $\bbu\in V_2(K)$ is approximated by
$$\bbu=\sum_{\sigma}\bbu_\sigma\varphi_\sigma+\xbar{\bbu}_K\;\xbar{\varphi}.$$  

%\begin{figure}[H]
%\begin{center}
%\begin{subfigure}{0.45\textwidth}
%\caption{$\bsigma$ is  a vertex}
%\includegraphics[trim=0.01cm 0.005cm 0.005cm 0.02cm,clip,width=4.5cm]{notationfigs/nplus}
%\end{subfigure}
%\begin{subfigure}{0.45\textwidth}
%\caption{$\bsigma$ is not a vertex}
%\includegraphics[trim=0.01cm 0.005cm 0.005cm 0.02cm,clip,width=4.5cm]{notationfigs/nmilieu}
%\end{subfigure}
%\end{center}
%\caption{\label{geometry} For any $\sigma\in K$, we define a {scaled} normal $\bbn_\sigma^K$ for the polygon $K$ and the DoFs $\sigma$ as follows: (a)--if $\sigma$ is a vertex, the vector $\bbn_\sigma^K$ is the sum of the {scaled outward} normals of the faces $e^+$ and $e^-$ sharing $\sigma$. (b)--If $\sigma$ is not a vertex, it belongs to the edge $e=e^+=e^-$ and $\bbn_\sigma^K$ is the outward normal to $e$.}
%\end{figure}

\begin{figure}[H]
\begin{center}
\subfigure[$\bsigma$ is  a vertex]{
\includegraphics[trim=0.01cm 0.005cm 0.005cm 0.02cm,clip,width=4.5cm]{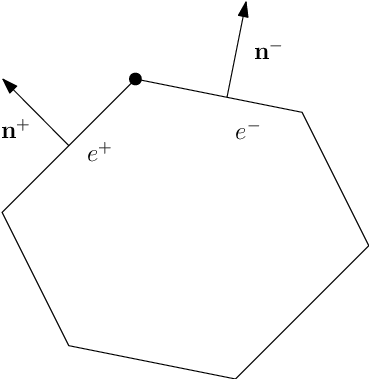}}
\subfigure[$\bsigma$ is not a vertex]
%\begin{subfigure}{0.45\textwidth}
{\includegraphics[trim=0.01cm 0.005cm 0.005cm 0.02cm,clip,width=4.5cm]{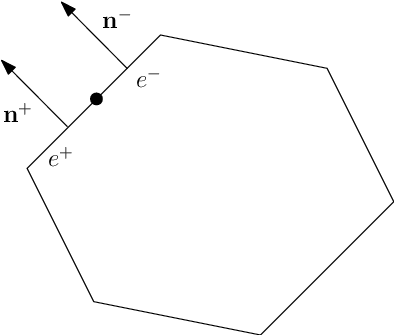}}
\end{center}
\caption{\label{geometry} For any $\sigma\in K$, we define a {scaled} normal $\bbn_\sigma^K$ for the polygon $K$ and the \dofs $\sigma$ as follows: (a)--if $\sigma$ is a vertex, the vector $\bbn_\sigma^K$ is the sum of the {scaled outward} normals of the faces $e^+$ and $e^-$ sharing $\sigma$. (b)--If $\sigma$ is not a vertex, it belongs to the edge $e=e^+=e^-$ and $\bbn_\sigma^K$ is the outward normal to $e$.}
\end{figure}
In \cite{BP_Pampa_VEM} the evolution of the degrees of freedom is done by:
\begin{subequations}\label{scheme}
\begin{enumerate}
\item Averaged value:
\begin{equation}\label{scheme:ave}\vert K\vert\,\dfrac{d\xbar{\bbu}_K}{dt}+\int_{\partial K} \bbf(\bbu)\cdot\bbn\; dl=0,\end{equation}
\item Point values
\begin{equation}\label{scheme:pt}\dfrac{d\bbu_{\sigma}}{dt}+\sum_{K, \sigma\in K}\bomega^K_{\sigma}\bbPhi_\sigma^K=0,\end{equation}
where
\begin{equation}
\label{scheme:pt:1}\bbPhi_\sigma^K=\bbJ_\sigma\cdot \big (\nabla \pi^0\bbu)({\bbx}_\sigma)+\bbPsi_\sigma^K
\end{equation}
with $\bbJ_\sigma$ the Jacobian of the flux evaluated at ${\bbx}_\sigma$ and $\bPsi_\sigma^K${ a filtering term we want to understand. We postpone its complete description until  section \ref{sec:analysis}. For now we only say its role is to remove some spurious wiggles  coming from the fact we may well have $\nabla\pi^0u=0$ with $u$ non constant on $K$.}

The matrix $\bomega_\sigma^K$ is defined as:
\begin{equation}\label{scheme:pt:2}\bomega_\sigma^K=\bbN_\sigma \bigg ( \big (\bbJ_\sigma\bbn_\sigma^K\big )^++\epsilon_\sigma^K\bigg)\qquad \epsilon_\sigma^K>0,
\end{equation}
$\bbn_\sigma^K$ is the outward normal to $K$ at $\sigma$ (see Figure \ref{geometry}) and 
\begin{equation}
\label{scheme:pt:3}\bbN_\sigma=\bigg ( \sum_{K, \sigma\in K} \bigg ( \big (\bbJ_\sigma\bbn_\sigma^K\big )^++\epsilon_\sigma^K\bigg)\bigg )^{-1}.
\end{equation}
If the system is equipped with an entropy, it is easy to show that if $\epsilon_\sigma^K>0$, this matrix is invertible, with a $L^2$ matrix norm uniformly bounded. We also see that
$$\sum_{K, \sigma\in K}\bomega_{\sigma}^K=\Id.$$
\end{enumerate}
\end{subequations}

The system \eqref{scheme} is evolved in time by a SSP Runge-Kutta scheme.
%%%%%%%%%%%%%%%%%%%
\subsection{Analysis of the \pampa scheme}\label{sec:analysis}
%%%%%%%%%%%%%%%%%%%

{Here, we recall some facts have been shown in \cite{pampa_dg,OEPampa} and demonstrate how to generalize them to polygons.} We will denote by $\langle \bullet  ,\bullet\rangle$ the  $L^2$ product
$$\langle u,v\rangle=\int_K u(\bbx)v(\bbx)\; d\bbx.$$
From Riesz representation theorem, we know that it is possible to find a family $\{\psi_\sigma\}\subset V_2(K)$ such that
$$\bbu(\sigma)=\langle \psi_\sigma, \bbu\rangle \text{ for all }\bbu\in V_2(K).$$
Since $\xbar \psi=1/\vert K\vert\in V_2(K)$ satisfies
$$\xbar{\bbu}_K=\langle \xbar \psi , \bbu\rangle,
$$ the family $\{ \{\psi_\bsigma\}_{\bsigma}, \xbar\psi\}$ is the  dual basis of $\{\{\varphi_\bsigma\},\xbar\varphi\}$ in $V_2(K)$. This means that we can recover the moments as $L^2$ inner products with functions in $V_2(K)$. Following \cite{pampa_dg,OEPampa}, we know that the dual functions $\{\psi_\sigma, {\xbar \psi}\}$ can be recovered from the basis function $\{\varphi_\bsigma , \bar \varphi\}$ from
\begin{equation}\label{dual basis}
M\begin{pmatrix} \vdots \\ \psi_\sigma\\ \vdots \\ \xbar \psi\end{pmatrix}=\begin{pmatrix} \vdots \\\varphi_\sigma\\ \vdots \\ \xbar \varphi\end{pmatrix}
\end{equation}
where $M$ is the mass matrix.
This is a positive definite matrix but unfortunately impossible to compute.

Let us write the \dg semi-discrete scheme on $K$. 
Since 
$$\frac{d\bbu_\sigma}{dt}=\langle \psi_\sigma , \frac{d\bbu}{dt}\rangle$$
\underline{\emph{if}} $\text{div }\bbf(\bbu)=\bbJ\cdot \nabla \bbu\in V_2(K)$, we can write
$$\frac{d\bbu_\sigma}{dt}=-{\frac{1}{\vert K\vert}}\int_K \psi_\sigma\text{ div }\bbf(\bbu)\; d\bbx=-\big (\bbJ\cdot \nabla \bbu\big )(\bbx_\sigma). $$ We also know that (whatever the flux)
$$\dfrac{d\xbar{\bbu}_K}{dt}=-{\frac{1}{\vert K\vert}}\int_{\partial K}\bbf(\bbu)\cdot\bbn\; dl=-{\frac{1}{\vert K\vert}}\int_K \text{ div }\bbf(\bbu)\; d\bbx,$$ so that \emph{if} $\text{div }\bbf(\bbu)\in V_2(K)$, 
and setting $\bbU{_K}=( \ldots , \bbu_\sigma, \ldots , \xbar{\bbu}_K)^T$, we could reinterpret the \pampa scheme as doing first
$$\begin{pmatrix}\vdots \\ S_K^\sigma\\\vdots \\ \xbar{S}_K\end{pmatrix}:=\dfrac{d\bbU_K}{dt}+M^{-1}\begin{pmatrix}\vdots\\
\int_K \varphi_\sigma\; \text{div }\bbf\; d\bbx\\\vdots \\
\int_K \xbar{\varphi}\; \text{div }\bbf\; dl
\end{pmatrix},$$
noticing that
$$\xbar{S}_K=\dfrac{d\xbar{\bbu}_K}{dt}+{\frac{1}{\vert K\vert}}\int_{\partial K} \bbf(\bbu)\cdot\bbn\; dl,$$
and then do
\begin{enumerate}
\item For the average value
$$\xbar{S}_K=\dfrac{d\xbar{\bbu}_K}{dt}+{\frac{1}{\vert K\vert}}\int_{\partial K} \bbf(\bbu)\cdot\bbn\; dl=0,$$
\item For the point values
$$\sum_{K, \sigma\in K}\bomega_\sigma^K S_\sigma^K=0.$$
\end{enumerate}
In other words, $S_\sigma^K=\dfrac{d\bbu_\sigma}{dt}+\bbPhi_\sigma^K$.

Unfortunately, in general $\bbJ\cdot \nabla \bbu\not\in V_2(K)$ because first $\nabla \bbu$ does not belong to $V_2(K)$. Hence we have this issue even for linear problems.  
Let us assume for a while that $\bbf(\bbu)=\bbA\bbu$ with $\bbA$ constant and denote by $\pi^0: V_k(K)\rightarrow\P_k(K)$ the $L^2$-orthogonal projector. We can write that

\begin{equation*}
\begin{split}
\int_K \psi_\sigma \bbA\cdot \nabla\bbu\; d\bbx&=\int_K\psi_\sigma \pi{^0}(\bbA\cdot\nabla\bbu)\; d\bbx+
\int_K \psi_\sigma \big (\Id-\pi{^0})\big ( \bbA\cdot\nabla\bbu\big )\; d\bbx\\&\quad=\int_K\pi{^0}(\psi_\sigma )\pi{^0}(\bbA\cdot\nabla\bbu)\; d\bbx+
\int_K(\Id-\pi{^0})(\psi_\sigma )\pi{^0}(\bbA\cdot\nabla\bbu)\; d\bbx\\
&\qquad +\int_K\big [ \pi{^0}\psi_\sigma \big ]\big [\big (\Id-\pi{^0})\big ( \bbA\cdot\nabla\bbu\big )\big ]\; d\bbx+\int_K\big [ \big (\Id-\pi{^0}\big )\psi_\sigma \big ]\big [\big (\Id-\pi{^0})\big ( \bbA\cdot\nabla\bbu\big )\big ]\; d\bbx.
\end{split}
\end{equation*}
Because $\pi{^{0}}$ is the  $L^2$ projection, we see that 
$$\int_K(\Id-\pi{^0})(\psi_\sigma )\pi{^0}(\bbA\cdot\nabla\bbu)\; d\bbx=0 \text{ and }\int_K\big [ \pi{^0}\psi_\sigma \big ]\big [\big (\Id-\pi{^0})\big ( \bbA\cdot\nabla\bbu\big )\big ]\; d\bbx=0$$
so that we have 
\begin{equation*}
\int_K \psi_\sigma \bbA\cdot \nabla\bbu\; d\bbx=\int_K\pi{^0}(\psi_\sigma )\pi{^0}(\bbA\nabla\bbu)\; d\bbx+
\int_K\big [ \big (\Id-\pi{^0}\big )\psi_\sigma \big ]\big [\big (\Id-\pi{^0})\big ( \bbA\nabla\bbu\big )\big ]\; d\bbx.
\end{equation*}
Defining $R^0=\Id-\pi^0$, we have
\begin{equation*}
 \int_K \psi_\sigma \bbA\cdot \nabla\bbu\; d\bbx=\int_K\pi{^0}(\psi_\sigma )\pi{^0}(\bbA\nabla\bbu)\; d\bbx+\int_KR{^0}(\psi_\sigma )R{^0}(\bbA\nabla\bbu)\; d\bbx
\end{equation*}
Because $\bbA$ is constant, we see that
$\pi{^0}(\bbA\nabla\bbu)=\bbA\cdot\pi{^0}(\nabla \bbu)=\bbA\cdot\nabla(\pi^0 \bbu)+O(h_K^3)$: Here we are making a $O(h^3_K)$ error since for all $\bbp\in {[\P_2(K)]^d}$,
$$\int_K \bbp\pi{^0}(\nabla u)\; d\bbx=\int_K\bbp \cdot\nabla u\; d\bbx=-\int_K\text{ div }\bbp \; u\; d\bbx+\int_{\partial K}\bbp\cdot\bbn \; u\;d\ell=\int_K \bbp\nabla \pi^0 u\; d\bbx+\int_{\partial K}\big (u-\pi^0 u)\bbp\cdot\bbn\; d\ell$$
since $\int_K (\pi^0-\Id)u\div\bbp \;d\bbx=0$ and $u-\pi^0 u=O(h_K^3)$.

Hence, we make the approximation
\begin{equation*}
\int_K \psi_\sigma \bbA\cdot \nabla\bbu\; d\bbx\remi{\approx }\big (\bbA \cdot \nabla\pi^0 \bbu\big )_\sigma+
\int_KR{^0}(\psi_\sigma) R{^0}( \bbA\nabla\bbu)\; d\bbx.
\end{equation*}

Following the VEM literature, we write
\begin{equation}\label{approx}
{\bbPsi_\sigma^K\approx
\alpha_K\vert K\vert \bigg ( \sum_{\sigma'} \big ((\Id-\pi^0)(\psi_\sigma)\big )_{\sigma'}\big ((\Id-\pi^0)(\nabla \bbu)\big )_{\sigma'}+\xbar{(\Id-\pi^0)\psi_\sigma}\; \xbar{(\Id-\pi^0)\nabla u}\bigg ).}\end{equation}
Then we make the additional approximation
$$\big ((\Id-\pi^0)(\psi_\sigma)\big )_{\sigma'}\big ((\Id-\pi^0)(\nabla \bbu)\big )_{\sigma'}\approx\frac{\big ((\Id-\pi^0)\psi_\sigma\big )_{\sigma'}\big ((\Id-\pi^0)( \bbu-\xbar{\bbu})\big )_{\sigma'}}{h_K}.$$
We note that $\xbar{ (\Id-\pi^0)\psi_\sigma}=0$ since
$$\xbar{\pi^0\psi_\sigma}=\frac{\int_K \pi^0\psi_\sigma \; d\bbx}{\vert K\vert }=
\frac{\int_K \psi_\sigma \; d\bbx}{\vert K\vert }=\xbar{\psi_\sigma}\;(=1)$$
because this is the value of the polynomials $1$ at $\sigma$.
Thus, there is no need to approximate the other term.
Last, in \eqref{approx},  $\alpha_K$ is an estimate of the maximum of the spectral radius of $\bbA\cdot\bbn$ for all unitary normals $\bbn$. Consequently, we have
\begin{equation}\label{eq:stabilized}
    \bbPsi_\sigma^K\approx\frac{\alpha_K\vert K\vert}{h_K}\sum_{\sigma'}\big ((\Id-\pi^0)\psi_\sigma\big )_{\sigma'}\big ((\Id-\pi^0)( \bbu-\xbar{\bbu})\big )_{\sigma'}
\end{equation}

Of course, it is not possible to compute the second term, because we have no direct access to the $\psi_\sigma$s, so we need to make a further approximation. There are two possibilities: either we approximate the $\psi_\sigma$ or we make an estimation:
\begin{itemize}
\item 
In the first case, to approximate the dual basis, we approximate the mass matrix. Let $\bbU_K$ be the vector of local VEM DoFs in the polygon $K$ and let
\begin{equation*}
    \Pi^0_K:\R^{N_k}\rightarrow\R^{N_p}
\end{equation*}
be the matrix that maps the VEM DoFs to the coefficients of $\pi^0\bbu$ in a polynomial basis $\{\mathcal{M}_\alpha\}_{\alpha=1}^{N_p}$. Define
\begin{equation*}
    (D_K)_{i\alpha}=\dof_i(\mathcal{M}_\alpha),\quad P_K^0=D_K\Pi_K^0,\quad R_K^0=\Id-P_K^0.
\end{equation*}
Let
\begin{equation*}
    (H_K)_{ij}=\int_K\mathcal{M}_i\mathcal{M}_j\;d\bbx
\end{equation*}
The computable VEM mass matrix is defined by
\begin{equation}\label{eq:VEM_mass}
    \M=(\Pi_K^0)^TH_K\Pi^0_K+\vert K\vert(R_K^0)^TR_K^0
\end{equation}
Then we evaluate $\M^{-1}$ and apply \eqref{dual basis} with $\M$ instead of $M$.
\item 
Another possibility, cheaper, is to note from \eqref{dual basis} that
$$\Vert\psi_\sigma\Vert_\infty\approx \frac{1}{\vert K \vert }$$
and then we replace $\big ((\Id-\pi^0)\psi_\sigma\big )(\bbx_\sigma')$ by $\frac{1}{\vert K\vert}\big ((\Id-\pi^0)\varphi_\sigma\big )(\bbx_\sigma')$, leading to 
\begin{equation}\label{approx:true}
{\bbPsi_\sigma^K=
\frac{\alpha_K}{h_K}\bigg ( \sum_{\sigma'} \big ((\Id-\pi^0)(\psi_\sigma)\big )_{\sigma'}\big ((\Id-\pi^0)(\bbu-\xbar{\bbu})\big )_{\sigma'}\bigg ).}\end{equation}
\end{itemize}
In practice, we use \eqref{approx:true} because it is simpler.

The scheme for point values is 
\begin{equation}
\label{scheme:ptbis}
\frac{d \bbu_\sigma}{dt}+\sum_{K, \sigma\in K} \bomega_{\sigma,K} \bigg ( \bbA\cdot \nabla \pi^0\bbu(\bbx_\sigma)+\bbPsi_\sigma^K\bigg )=0,
\end{equation}
and nothing is to be changed for the average value due to the structure of the last line of $M^{-1}$: since $1=\sum_\bsigma \varphi_\bsigma+\xbar\varphi$, we see that the last line of $M^{-1}$ must be made of $1$. Using the symmetry of $M$, the inverse has the structure:
$$M^{-1}=\begin{pmatrix} B & \mathbf{1}^T\\
\mathbf{1}& 1\end{pmatrix}, \quad B^T=B.$$

%%%%%%%%%%%%%%%%%%%
\section{Boundary conditions}
%%%%%%%%%%%%%%%%%%%
For the boundary conditions, again we make use of the same arguments.
Writing first that 
$$\int_K \psi_\bsigma \text{div }\bbf(\bbu)\; d\bbx=-\int_K\nabla\psi_\bsigma \cdot \bbf(\bbu)\; d\bbx+\int_{\partial K}\psi_\bsigma \bbf(\bbu)\bbn\; d\gamma$$
if we want to impose a boundary condition by flux on one of the faces $e_b$ of $\partial K$ (i.e. $e_b\subset K\cap\partial K$), we simply have to add 
$$\int_{e_b}\psi_\bsigma \big ( \hbbf_\bbn(\bbu, \bbu_b)-\bbf(\bbu)\bbn\big )\; d\gamma$$ to the above term. Again we face the same problems of computability. One way to proceed is the following:
Since $\varphi_{\xi}$ vanishes on $e_b$ if ${\xi\not\in e_b}$, we introduce the vector 
$$\bbX=\bigg ( \ldots \int_{e_b}\varphi_\xi \big ( \hbbf_\bbn(\bbu, \bbu_b)-\bbf(\bbu)\bbn\big )\; d\gamma, \ldots, \int_{e_b}\xbar{\varphi} \big ( \hbbf_\bbn(\bbu, \bbu_b)-\bbf(\bbu)\bbn\big )\; d\gamma.$$
All the terms that correspond to \dofs that are not on $e_b$ are equal to $0$.
If  we could evaluate $M^{-1}X$, we would get the vector of components
$$\bigg (\ldots \int_{e_b}\psi_\bsigma \big ( \hbbf_\bbn(\bbu, \bbu_b)-\bbf(\bbu)\bbn\big )\; d\gamma,\ldots, \int_{e_b}\big ( \hbbf_\bbn(\bbu, \bbu_b)-\bbf(\bbu)\bbn\; d\gamma\bigg )$$
Since $M^{-1}$ is impossible to evaluate, {and we  replace it by $\M^{-1}$ that defined in \eqref{eq:VEM_mass}.}

Hence, 
 we take
$$\bbPhi_\bsigma^K=\bbJ(\bbu_\bsigma)(\nabla\pi^0\bbu)(\bbx_\bsigma)+\bbPsi_\sigma^K+\int_{e_b}\widetilde{\psi_\bsigma}\big ( \hbbf_\bbn(\bbu, \bbu_b)-\bbf(\bbu)\bbn\big )\; d\gamma,$$
where $\widetilde{\psi_\bsigma}$ is  the approximated dual function $\psi_\bsigma$.
\begin{remark}[Reduction to the case of meshes with simplex]
Of course, when we consider meshes with simplex and use the basis function, we can apply the same technique as here with the exact knowledge of the mass matrix. The stabilisation term $\bbPsi_\bsigma^K$ is not needed.
\end{remark}

%%%%%%%%%%%%%%%%%%%
\section{Nonlinear stabilization for \pampa-VEM}
%%%%%%%%%%%%%%%%%%%
It is known that the standard third-order \pampa (semi-discrete Active Flux) method will not be oscillation free without any form of limiting strategies. When the solution
develops discontinuities, these schemes will be prone to numerical oscillations and possibly non-physical solutions. For instance, in
gas dynamics, the density or the internal energy, or both might become negative. To address this issue, limiting is required. In \cite{AbgrallBoscheriLiu2025}, we have introduced the first-order PAMPA scheme and a convex monolithic limiting approach to blend the high- and low-order schemes so that positivity property is ensured. 

Inspired by the oscillation-free and oscillation elimination DG methods \cite{Lu2021,Peng2025}, we have introduced an effective convex oscillation elimination (COE) parameter to suppress spurious oscillation near shocks and strong discontinuities in \cite{OEPampa} on unstructured triangle meshes. In this section, we show how to generalize this parameter to the polygonal meshes. The element-based COE parameter is given by
\begin{equation*}
  \theta_K=\exp\Big(-\sum_{e\subset\partial K}\frac{\alpha_e\dt}{h_{e,K}}\sigma_{e,K}(\bbu)\Big),\quad h_{e,K}=\sup_{\bbx\in K}\text{dist}~(\bbx, e),\quad \delta_{e,K}=\max_{i}\sigma_{e,K}(\bbu^{(i)})
\end{equation*}
\begin{equation*}
  \sigma_{e,K}(\bbu^{(i)})=\left\{\begin{aligned}
&0, &&\mbox{if}~\bbu^{(i)}\equiv\langle\bbu^{(i)}\rangle_\Omega\\
&\sum_{\kappa=1}^{k=2}h_{e,K}^\kappa\sum_{\vert\bm{\alpha}\vert=\kappa}\dfrac{\dfrac{1}{\vert e\vert}\int_e\vert \llbracket \partial^{\bm{\alpha}}\big(\pi^0\bbu^{(i)}\big)\rrbracket_e\vert\;{\mathrm d}\ell}{\Vert \big(\pi^0\bbu^{(i)}\big)-\langle \big(\pi^0\bbu^{(i)}\big)\rangle_\Omega\Vert_{L^\infty(\Omega)}},&&\mbox{otherwise.}
\end{aligned}\right.
\end{equation*}

\noindent Last, in \cite{abgrall2026constructionentropysatisfyingactive}, we have shown how to construct semi-discrete entropy satisfying scheme.

Each of these strategies translate by blending a first order scheme and a high order scheme (the one described above). The first order scheme is of Local Lax-Friedrichs type. One can combine each approach, and we get a Bound Preserving, Oscilation Eliminated and Entropy satisfying scheme by using a worst case scenario. This is how we have obtained the numerical results of the next section.
%%%%%%%%%%%%%%%%%%%
\section{Numerical experiments}
%%%%%%%%%%%%%%%%%%%
\subsection{Shock Vortex interaction}

%\begin{figure}[ht!]
%\centerline{\subfigure[$\rho(\bbx,t=)$]{\includegraphics[trim=4.2cm 2.5cm 1.5cm 18.0cm,clip,width=5.25cm]{Euler/ShockVortex/ro1.png}}\hspace*{0.05cm}
%	\subfigure[$\rho(\bbx,t=)$]{\includegraphics[trim=4.2cm 2.5cm 1.5cm 18.0cm,clip,width=5.25cm]{Euler/ShockVortex/ro5.png}}\hspace*{0.05cm}
%	\subfigure[$\rho(\bbx,t=)$]{\includegraphics[trim=4.2cm 2.5cm 1.5cm 18.0cm,clip,width=5.25cm]{Euler/ShockVortex/ro9.png}}}
%    \vskip1pt
%    \centerline{\subfigure[$t=0.529$]{\includegraphics[trim=4.2cm 2.5cm 1.5cm 18.0cm,clip,width=5.25cm]{Euler/ShockVortex/sens1.png}}\hspace*{0.05cm}
%	\subfigure[$t=0.662$]{\includegraphics[trim=4.2cm 2.5cm 1.5cm 18.0cm,clip,width=5.25cm]{Euler/ShockVortex/sens5.png}}\hspace*{0.05cm}
%	\subfigure[$t=0.8$]{\includegraphics[trim=4.2cm 2.5cm 1.5cm 18.0cm,clip,width=5.25cm]{Euler/ShockVortex/sens9.png}}}
%\caption{\sf Example : Top: isolines of the density at 3 instances of the simulation.
%We have plotted 30 equaly spaces isolines, $\rho\in [0.88,1.181]$.
%Bottom: plot of the sensor. Red corresponds to $1$, and yellow is about $0.75$. \label{fig: Shock_Vortex}}
%\end{figure}
%
The domain is $[0,2]\times [0.2,0.8]$. The left are right faces are inflow/outflow, the top and bottom ones are walls. The mesh is filled with quad seen as polygons. It contains $323,220$ \dofs and $107,221$ elements. This corresponds to a $180\times 600$ mesh.

\begin{figure}[ht!]
\centerline{\subfigure[$\rho(\bbx,t=0.529)$]{\includegraphics[trim=4.2cm 2.5cm 1.5cm 18.0cm,clip,width=5.25cm]{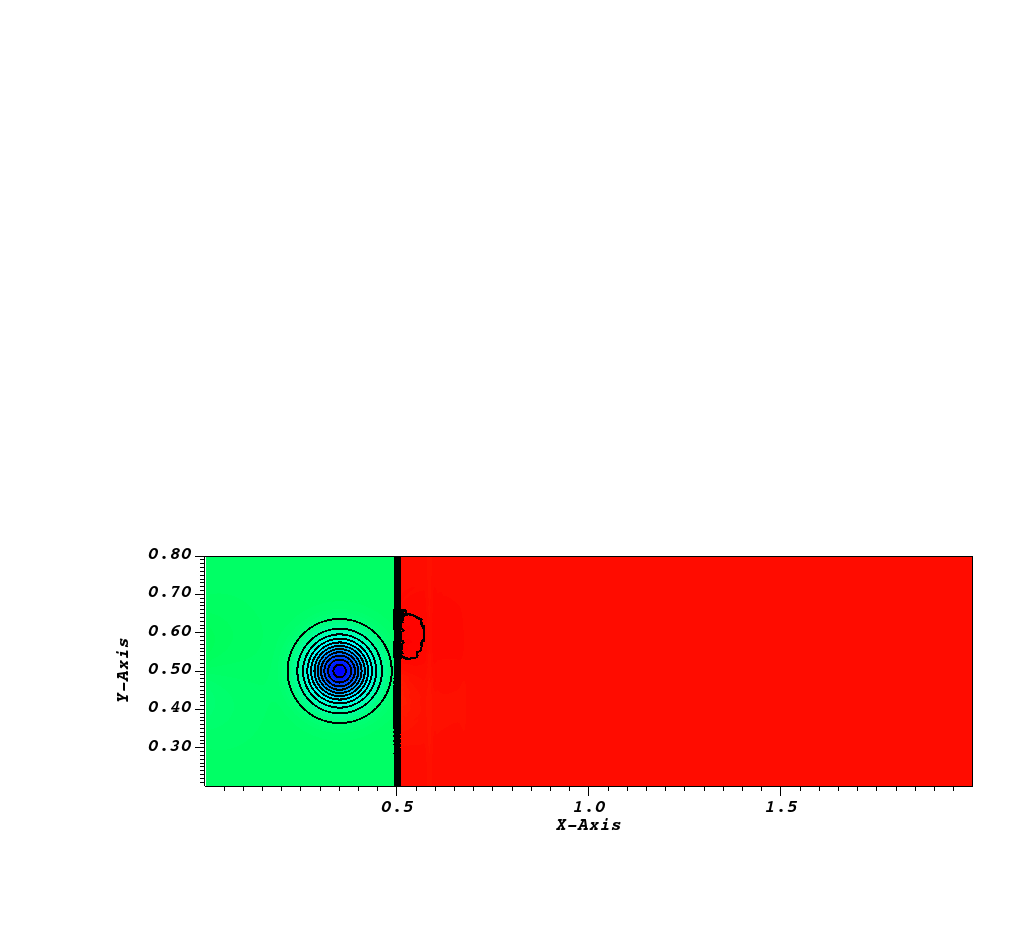}}\hspace*{0.05cm}
	\subfigure[$\rho(\bbx,t=0.662)$]{\includegraphics[trim=4.2cm 2.5cm 1.5cm 18.0cm,clip,width=5.25cm]{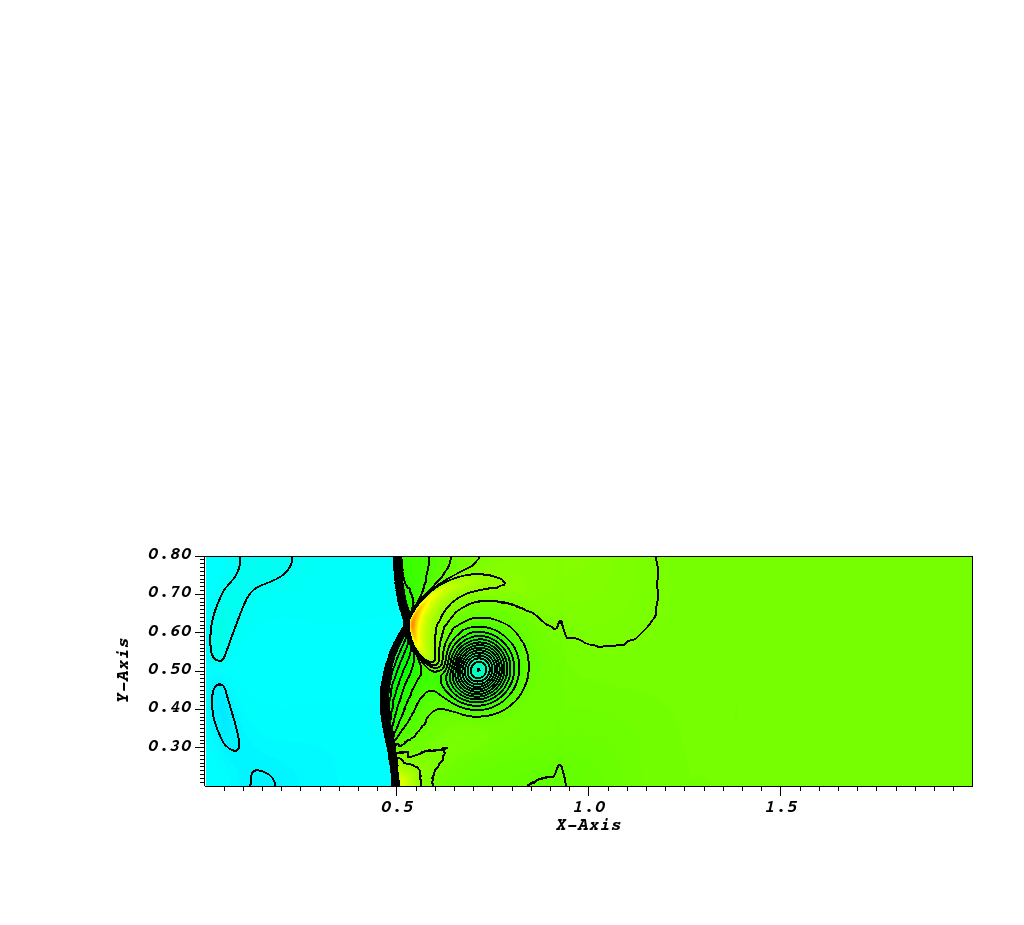}}\hspace*{0.05cm}
	\subfigure[$\rho(\bbx,t=0.8)$]{\includegraphics[trim=4.2cm 2.5cm 1.5cm 18.0cm,clip,width=5.25cm]{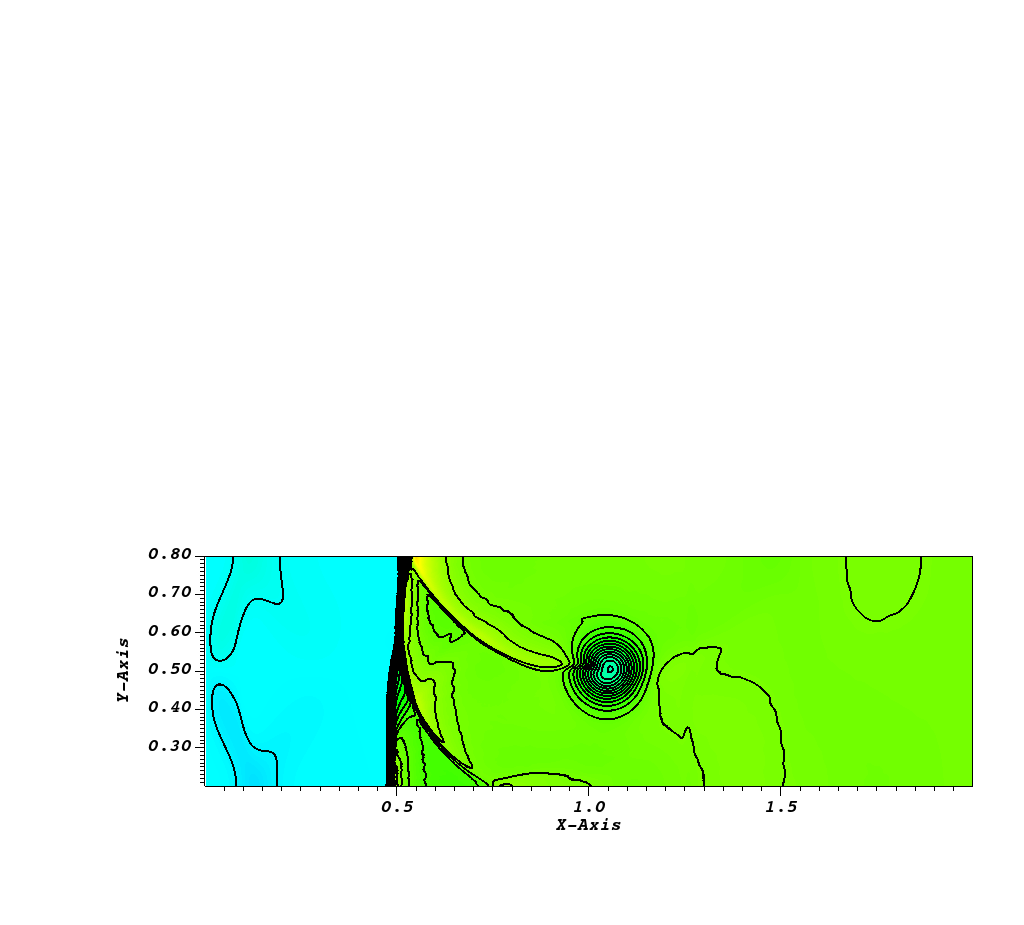}}}
    \vskip1pt
    \centerline{\subfigure[$t=0.529$]{\includegraphics[trim=4.2cm 2.5cm 1.5cm 18.0cm,clip,width=5.25cm]{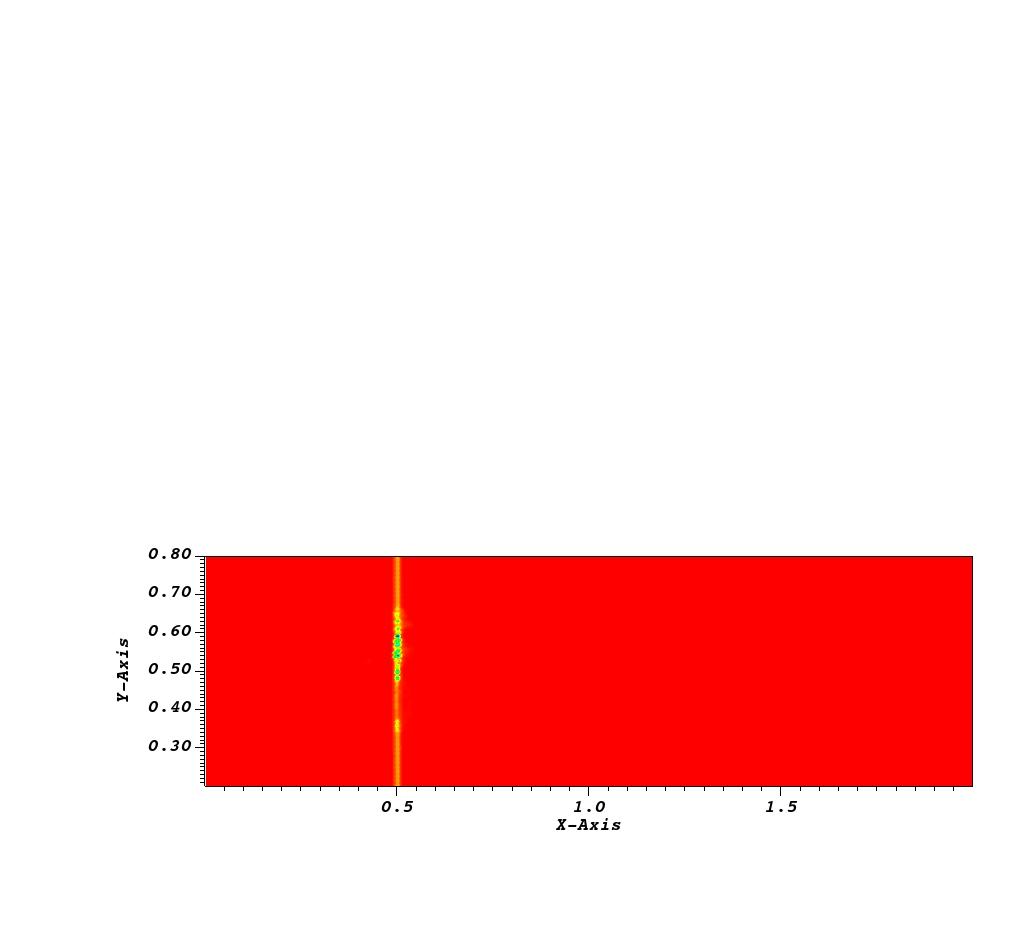}}\hspace*{0.05cm}
	\subfigure[$t=0.662$]{\includegraphics[trim=4.2cm 2.5cm 1.5cm 18.0cm,clip,width=5.25cm]{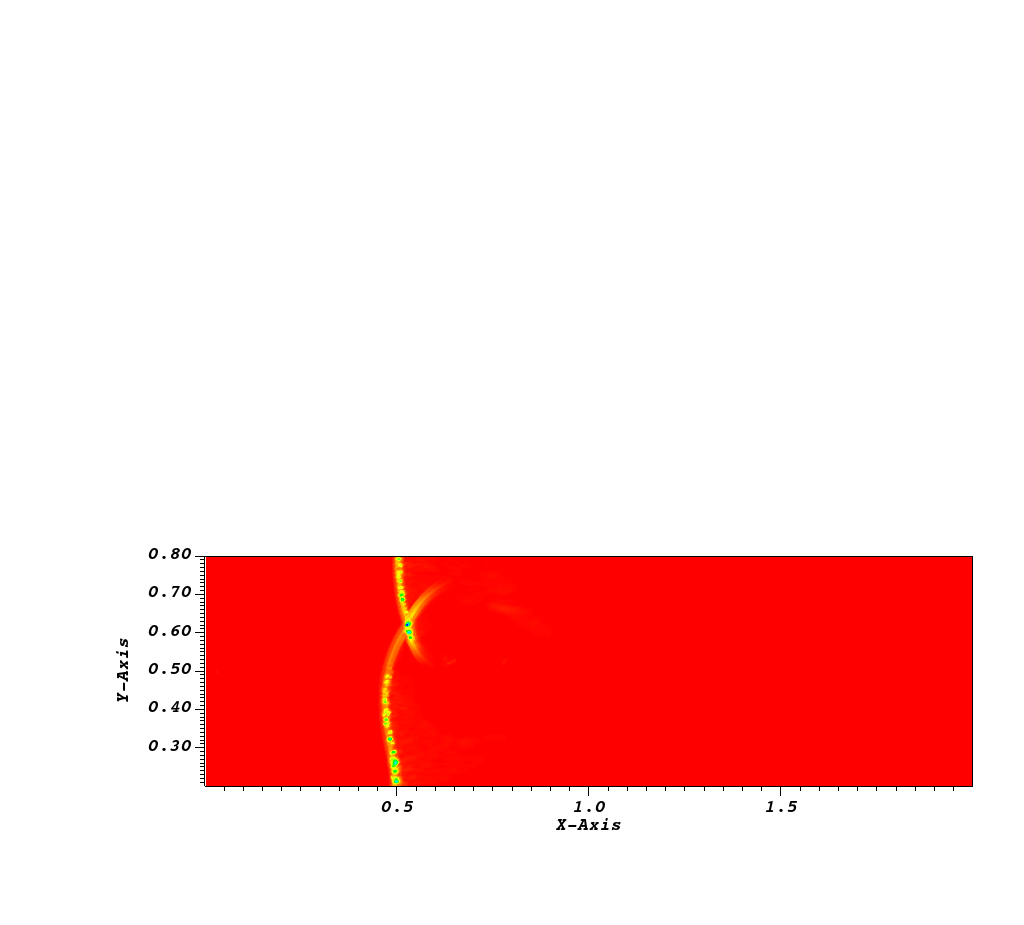}}\hspace*{0.05cm}
	\subfigure[$t=0.8$]{\includegraphics[trim=4.2cm 2.5cm 1.5cm 18.0cm,clip,width=5.25cm]{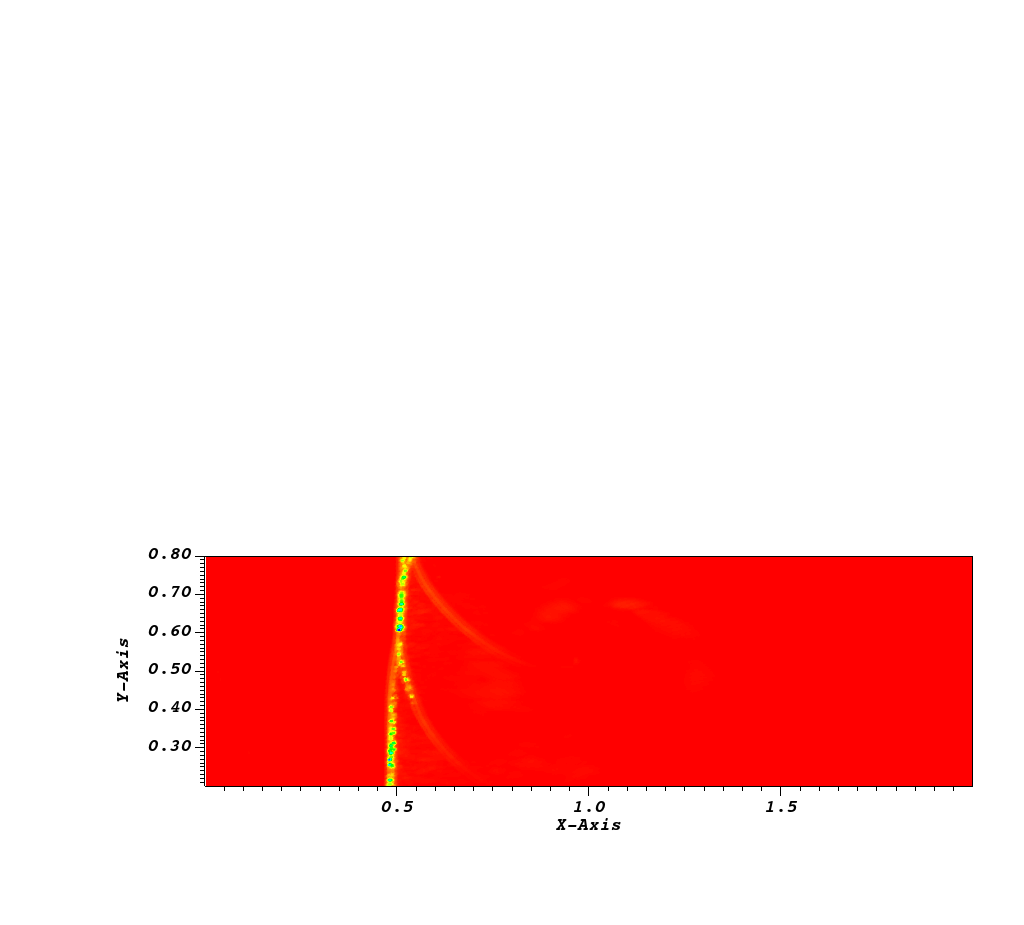}}}
\caption{\sf Example : Top: isolines of the density at 3 instances of the simulation.
We have plotted 30 equaly spaces isolines, $\rho\in [0.88,1.181]$.
Bottom: plot of the sensor. Red corresponds to $1$, and yellow is about $0.75$. \label{fig: Shock_Vortex}}
\end{figure}

What can be observed is that a fraction of the first order is only applied in and localy around the shocks. The vortex is not affected by this. This is the challenging aspect of this simulation.

\subsection{Shock diffraction over a wedge}
\begin{figure}[ht]
\begin{center}
\subfigure[]{
\includegraphics[width=0.45\textwidth]{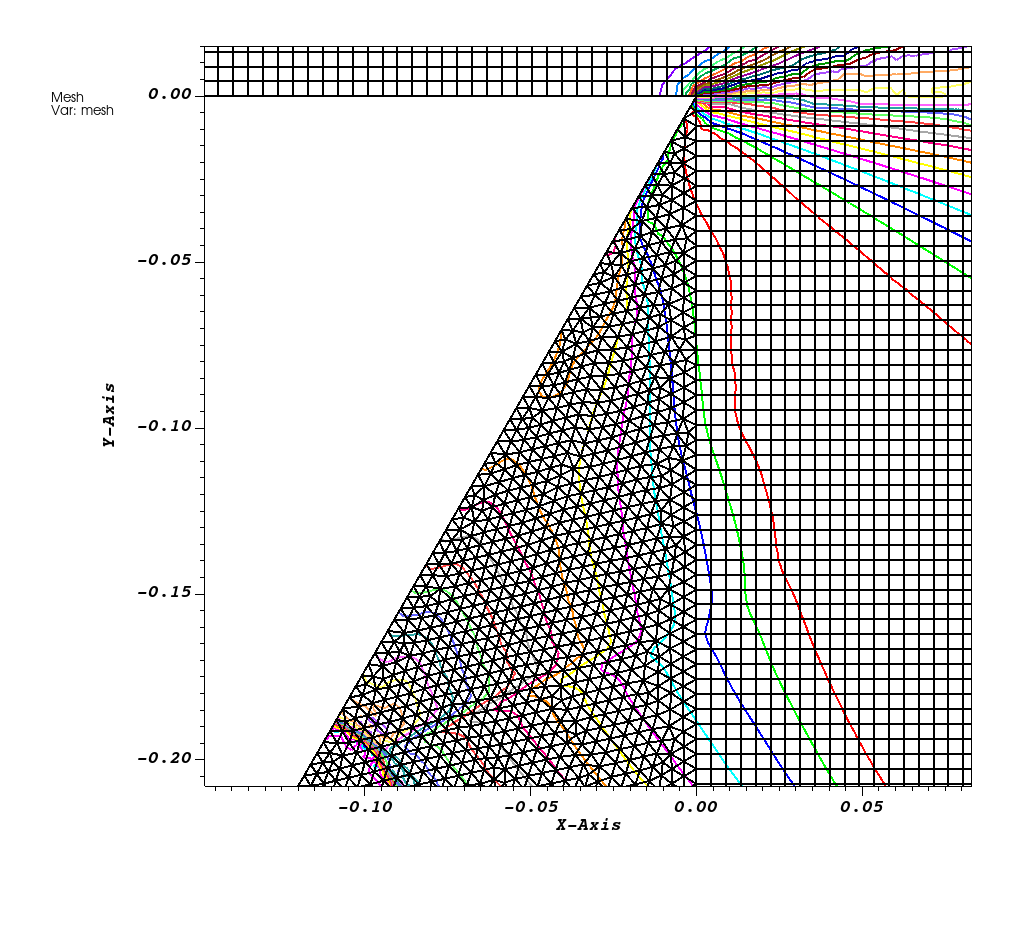}}
\subfigure[]{\includegraphics[width=0.45\textwidth]{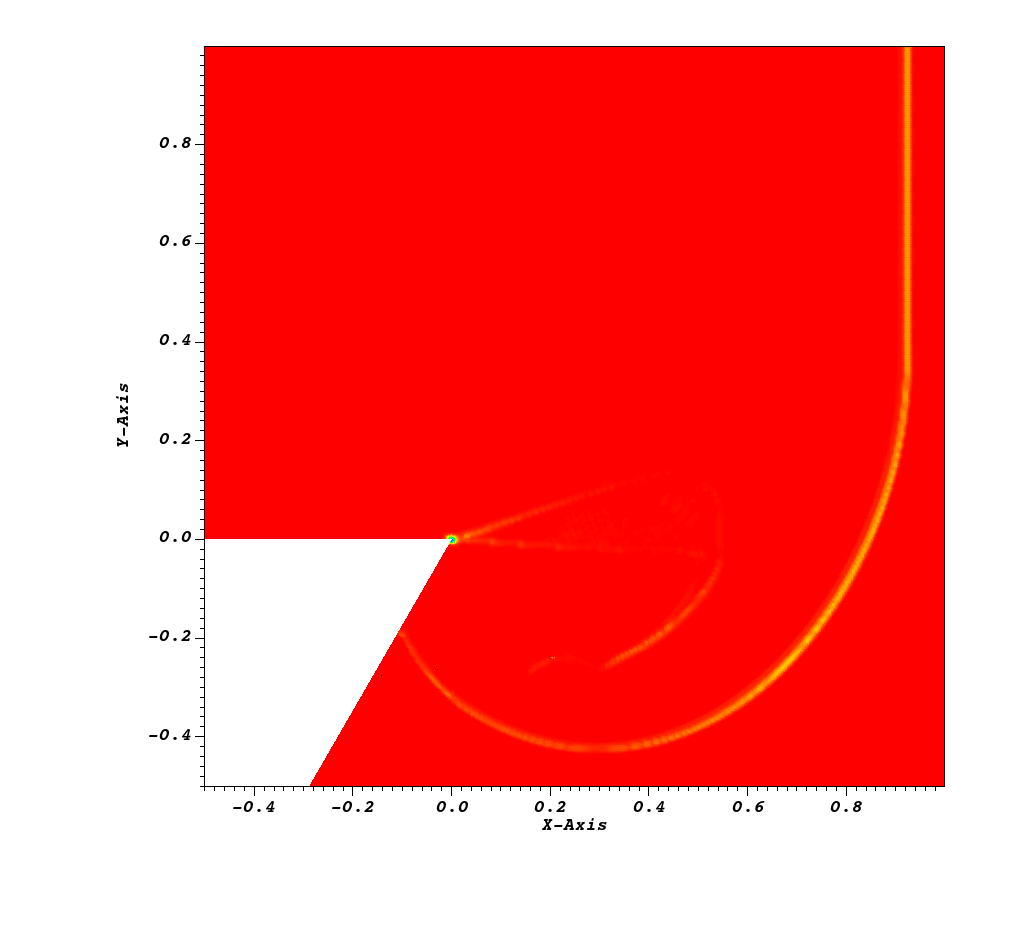}}
\end{center}
\caption{\label{fig:diff:mesh}(a): Zoom of the mesh for the shock diffraction problem. The isolines of the density, as in figure \ref{fig:shock_diff} are shown.(b): plot of the sensor, blue is 0, red is 1.}
\end{figure}

\begin{figure}[ht!]
\centerline{\subfigure[$\rho(\bbx,t=8\, 10^{-2})$]{\includegraphics[trim=4.2cm 2.5cm 1.5cm 0.01cm,clip,width=0.4\textwidth]{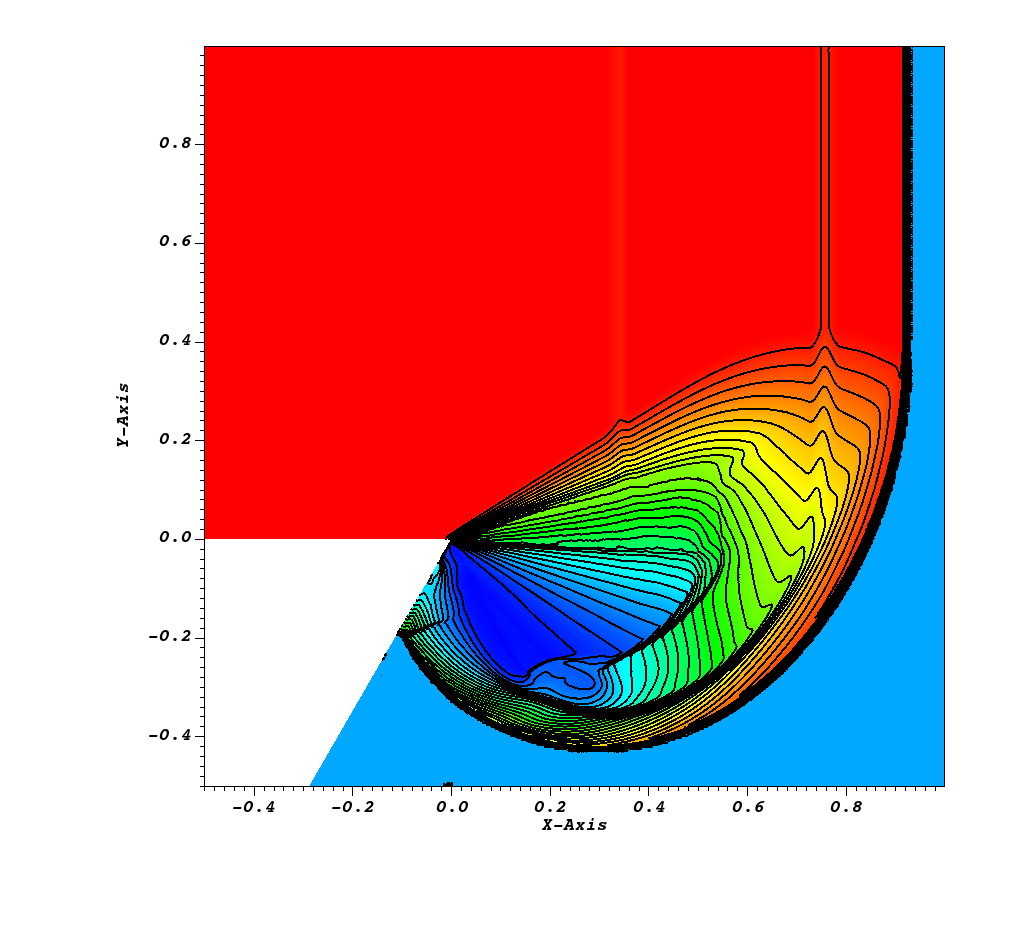}}\hspace*{1.5cm}
	\subfigure[$p(\bbx,t=8\, 10^{-2})$]{\includegraphics[trim=4.2cm 2.5cm 1.5cm 0.01cm,clip,width=0.4\textwidth]{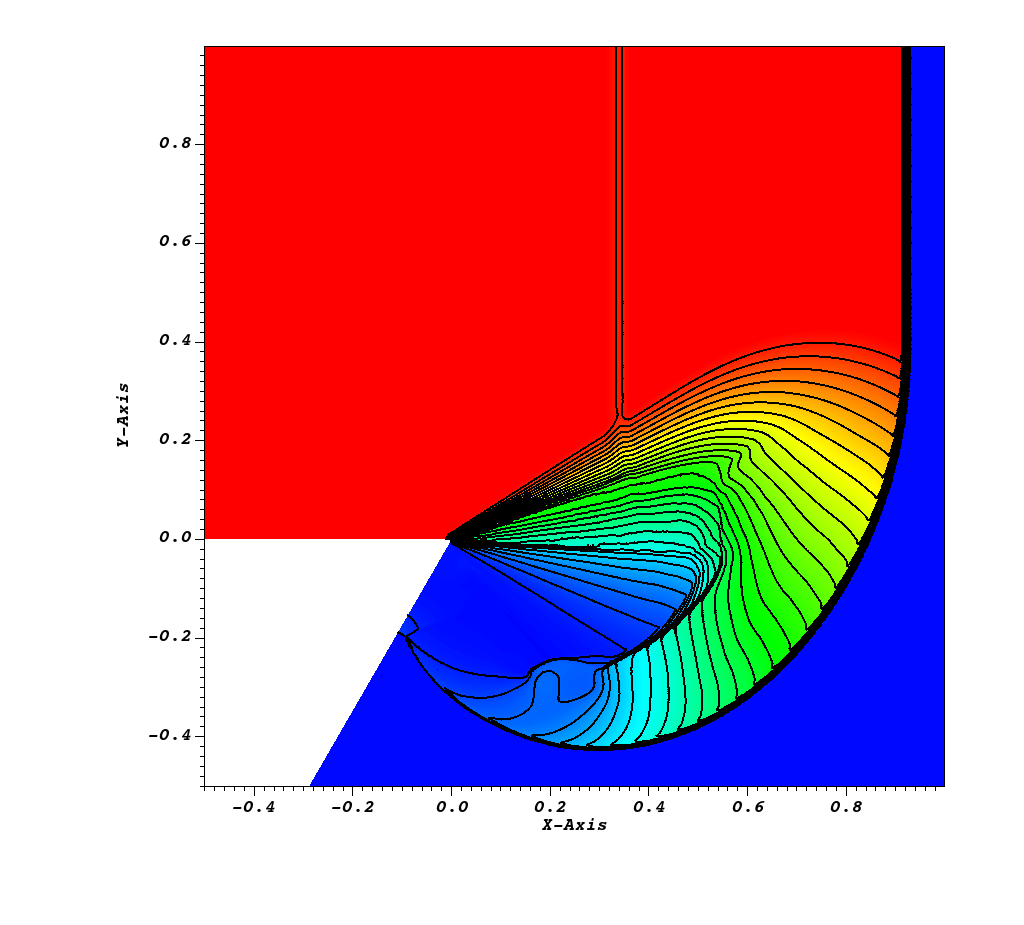}}}
\caption{\label{fig:shock_diff}Shock Diffraction problem, (a): Density 30 equally spaced isolines, $\rho\in [0.07584,5.735]$, (b): Pressure, 30 equally spaced isolines,  $p\in [0.03788,116.9]$}
\end{figure}

The problem is a shock diffraction over a wedge with angle 120$^\circ$. The domain is $[-0.5,1]\times [0,1]\cup[0,1]\times [0,-0.5]\cup T$ where $T$ is the triangle of vertices $(0,0)$, $(0,-0.5)$, $-\frac{\sqrt{3}}{6},-0.5)$.  The lect and right faces of the domain are inflow/outflow, the other boundary are all walls. On the right, the state is $(\rho, u,v,p)=(1,0,0,\gamma)$ with $\gamma=1.4$. The left state is defined by the Rankine Hugoniot conditions corresponding to a $\text{Mach}=10$ shock. The initial shock is localised at $x=-0.025$
 The mesh is represented  in figure \ref{fig:diff:mesh}, and the solution on figure \ref{fig:shock_diff}. It  contains quadrangles and triangles, all are seen as polygons. There are 108,280 polygons and 317,887 \dofs. The simulation is run until $T=8.\, 10^{-3}$. 

Again, one can observe that a fraction of the first order scheme is only applied near the corner and the strongest shocks. The rest of solution is virtualy high order.

\section{Conclusions and perspectives}
We have clarified the way to implement boundary conditions and the need for a filtering mechanism when considering the semi-discrete Active flux type scheme \pampa on a polygonal mesh.
The construction can be easily extended to higher order approximation following exactly the same procedure, simply because the internal moments that are needed are naturally represented as the $L^2$ moments of some polynomial function, exactly as the average is represented as the moment with respect to the constant function $1$.

 \printbibliography
\end{document}